\documentclass[a4paper,oneside,british]{amsart}
\usepackage[T1]{fontenc}
\usepackage[latin1]{inputenc}
\usepackage{amssymb}

\makeatletter
 \theoremstyle{plain}    
 \newtheorem*{prop*}{Proposition} 
 \theoremstyle{plain}    
 \newtheorem{lem}{Lemma} 
 \theoremstyle{remark}
 \newtheorem*{rem*}{Remark}
 \theoremstyle{plain}    
 \newtheorem{prop}{Proposition} 
 \theoremstyle{plain}    
 \newtheorem{thm}{Theorem} 

\usepackage{alex}

\usepackage{babel}
\makeatother
\begin{document}
\title[Representations and Extended Deligne-Lusztig Varieties]{Representations of Reductive Groups over Finite Rings and Extended Deligne-Lusztig Varieties}

\author{Alexander Stasinski}

\address{School of Mathematical Sciences\\
University of Nottingham\\
NG7 2RD\\
England}

\email{pmxas@maths.nott.ac.uk}

\begin{abstract}
In a previous paper it was shown that a certain family of varieties
suggested by Lusztig, is not enough to construct all irreducible complex
representations of reductive groups over finite rings coming from
the ring of integers in a local field, modulo a power of the maximal
ideal. In this paper we define a generalisation of Lusztig's varieties,
corresponding to an extension of the maximal unramified extension
of the local field. We show in a particular case that all irreducible
representations appear in the cohomology of some extended variety.
We conclude with a discussion about reformulation of Lusztig's conjecture.
\end{abstract}
\maketitle

\section{Introduction}

Let $F$ be a local field with finite residue field $\mathbb{F}_{q}$,
ring of integers $\mathcal{O}_{F}$, and maximal ideal $\mathfrak{p}_{F}$.
Let $\mathbf{G}$ be a connected reductive group over $F$, and set
$G_{r}=\mathbf{G}(\mathcal{O}_{F}/\mathfrak{p}_{F}^{r})$, for any
integer $r\geq1$.

Let $F^{\textrm{ur}}$ be a maximal unramified extension of $F$,
and let $\mathbf{G}_{r}=\mathbf{G}(\mathcal{O}_{F^{\textrm{ur}}}/\mathfrak{p}_{F^{\textrm{ur}}}^{r})$.
Then $\mathbf{G}_{r}$ carries a structure of affine algebraic group
over the residue field $\mathbb{F}$ of $F^{\textrm{ur}}$, and is
equipped with a morphism $\phi:\mathbf{G}_{r}\rightarrow\mathbf{G}_{r}$
induced by the Frobenius automorphism in $\textrm{Gal}(F^{\textrm{ur}}/F$),
such that $\mathbf{G}_{r}^{\phi}=G_{r}$.

In a previous version of \cite{lusztig-preprint-published}, Lusztig
conjectured (in the setting where $F$ is of equal characteristic
and $\mathbf{G}$ defined over $\mathbb{F}_{q}$) that every irreducible
representation of $G_{r}$ appears in the $l$-adic cohomology of
a variety\[
X_{x}=\{ g\in\mathbf{G}_{r}\mid g^{-1}\phi(g)\in x\mathbf{U}_{r}\},\]
for some element $x\in\mathbf{G}_{r}$. The results in \cite{alex_paper1}
show that this does not hold for the case $\mathbf{G}=\textrm{SL}_{2}$,
$F$ of equal characteristic, $r=2$, and $q$ odd, but that the representations
not accounted for by Lusztig's varieties, can in this case all be
realised by a variety of a different kind. The latter variety was
however constructed in a rather ad hoc manner, and it did not make
it clear how to generalise it. It is therefore natural to seek a construction
that extends Lusztig's conjecture in a more uniform and conceptual
way. We believe that the ideas in this paper provide a first step
in this direction.

The following is an outline of the contents. First we review Greenberg's
theory of reduction of schemes over local rings, modulo some power
of the maximal ideal. This provides the general theory underlying
the construction and properties of the groups $\mathbf{G}_{r}$, and
their generalisations.

Next, we define groups $\mathbf{G}_{L,r}$ for any finite Galois extension
$L/F^{\textrm{ur}}$. The fact that $G_{r}$ can be realised as the
fixed points of $\mathbf{G}_{L,r}$ under the actions of the elements
of the Galois group $\textrm{Gal}(L/F)$, allows us to define a generalisation
of the varieties of Lusztig, which we call extended Deligne-Lusztig
varieties.

The main result of this paper is that every irreducible representation
of dimension $(q^{2}-1)/2$ of the group $G_{2}$, where $\mathbf{G}=\textrm{SL}_{2}$,
$F$ of equal characteristic, $r=2$, and $q$ odd, is realised by
a certain extended Deligne-Lusztig variety, which is an analogue of
Lusztig's variety $X_{x}$ for $x=\begin{pmatrix}1 & 0\\
\epsilon & 1\end{pmatrix}$, where $\epsilon$ is a prime element in $\mathcal{O}_{F^{\textrm{ur}}}$.
We show that the ad hoc variety constructed in \cite{alex_paper1}
is in fact isomorphic to the quotient of this extended variety, modulo
a finite group. Together with the results of Lusztig (cf. \cite{lusztig-preprint-published},
sect. 3) this shows that all irreducible representations of $G_{2}$
occur in the cohomology of some extended Deligne-Lusztig variety of
a certain kind.

In the final section we speculate on how the construction of extended
Deligne-Lusztig varieties may be used to reformulate Lusztig's conjecture.
This perspective is still very rudimentary, which is shown by the
large number of open questions.

This work will form part of the author's PhD thesis, and has been
carried out at the University of Nottingham under the supervision
of Prof. Ivan Fesenko.

\section{Reduction of schemes over local rings}

In order to study the representation theory of reductive groups over
the ring of integers in a local field modulo some power of the maximal
ideal, it is very useful (perhaps essential) to view these groups
as fixed point subgroups of certain algebraic groups which are reductions
of the reductive group in question, modulo a power of the maximal
ideal. Even though this theory of reduction can be described in elementary
terms for affine varieties, it is more convenient to state it in the
general case of schemes of finite type over certain local rings. The
theory in this form, and all of the results in this section, are due
to Greenberg (cf. \cite{greenberg1}, \cite{greenberg2}).

Let $R$ be a commutative noetherian local ring, with maximal ideal
$\mathfrak{m}$ and residue field $k$, assumed to be perfect in the
mixed characteristic case. Denote by $X$ a scheme of finite type
over $R$. Then there exists a functor $\mathcal{F}_{R}$ from the
category of schemes of finite type over $R$, to schemes of over $k$,
and a functor $\mathcal{G}_{R}$ from schemes over $k$ to schemes
over $R$, which is adjoint to $\mathcal{F}_{R}$ in the sense that\[
\Hom_{R}(\mathcal{G}_{R}(Y),X)\cong\Hom_{k}(Y,\mathcal{F}_{R}(X)),\]
for any scheme $Y$ over $k$.

The functor $\mathcal{G}_{R}$ has the property that $\mathcal{G}_{R}(\Spec(k))=\Spec(R)$,
and so $X(R)\cong(\mathcal{F}_{R}X)(k)$, i.e. there exists a bijection
between the points of $X$ with values in $R$, and the points of
$\mathcal{F}_{R}(X)$ with values in $k$.

An advantage of the functorial approach is that it is straightforward,
once the categorical framework has been introduced, to show that the
reduction of a group scheme is itself a group scheme. More precisely:

\begin{prop*}
If $X$ is a group scheme over $R$, then $\mathcal{F}_{R}X$ is a
group scheme over $k$, and for every scheme $Y$ over $k$, the bijection
$(\mathcal{F}_{R}X)(Y)\cong(\mathcal{G}_{R}Y)(X)$ is an isomorphism
of groups.
\end{prop*}
\begin{proof}
This is Corollary 1, p. 639, and Corollary 4, p. 641 in \cite{greenberg1}.
\end{proof}
Furthermore, the functor $\mathcal{F}_{R}$ preserves subschemes,
takes schemes of finite type to schemes of finite type, separated
schemes to separated schemes, and affine schemes to affine schemes,
(cf. \cite{greenberg1}, Theorem, \S4).

Now suppose we have a second commutative noetherian local ring $R'$,
and a homomorphism $\phi:R\rightarrow R'$. Let $X^{\phi}$ be the
scheme over $R'$ obtained from $X$ by extension of scalars. Then
another consequence of the functorial construction is the existence
of a connecting morphism\[
\mathcal{F}_{\phi}(X):\mathcal{F}_{R}X\longrightarrow\mathcal{F}_{R'}X^{\phi},\]
which is a group homomorphism if $X$ is a group scheme, (cf. \cite{greenberg1},
\S5).

Assume from now on that $\phi$ is surjective, and that the ideal
$\Ker\phi$ is annihilated by $\mathfrak{m}$. Let $\rho:R\rightarrow R/\mathfrak{m}$
be the canonical map. We call a scheme \emph{simple} over $R$ \emph{}if
it has non-singular non-degenerate reduction mod $\mathfrak{m}$\emph{.}
Then we have the following results:

\begin{prop*}
If $X$ is simple over $R$, then $\mathcal{F}_{\phi}(X)$ is surjective.
\end{prop*}
\begin{proof}
Cf. Corollary 2, p. 262 in \cite{greenberg2}.
\end{proof}
\begin{prop*}
If $X$ is simple over $R$ and $X^{\rho}$ is reduced and irreducible,
then $\mathcal{F}_{R}X$ is reduced and irreducible.
\end{prop*}
\begin{proof}
Cf. Corollary 2, p. 264 in \cite{greenberg2}.
\end{proof}
We will now say something about projective limits. Suppose that $R$
is complete with respect to its $\mathfrak{m}$-adic topology, i.e.
$R\cong\prolim R/\mathfrak{m}^{i}$, with the connecting morphisms
$\rho_{i,j}:R/\mathfrak{m}^{i}\rightarrow R/\mathfrak{m}^{j}$, for
$i\geq j$. Let $X_{\mathfrak{m}^{i}}$ be the scheme over $R/\mathfrak{m}^{i}$
obtained from $X$ by extension of scalars, and set $\mathcal{F}_{\mathfrak{m}^{i}}X:=\mathcal{F}_{R/\mathfrak{m}^{i}}X_{\mathfrak{m}^{i}}$.
Then it is shown in \cite{greenberg1}, \S6 that the induced maps
$\mathcal{F}_{\rho_{i,j}}(X)$ form a projective system, and there
is a bijection\[
X(R)\longrightarrow\prolim(\mathcal{F}_{\mathfrak{m}^{i}}X)(k),\]
which is functorial in $X$. Moreover, if $X$ is a group scheme,
so is each $\mathcal{F}_{\mathfrak{m}^{i}}X$, each map $\mathcal{F}_{\rho_{i,j}}(X)$
is a homomorphism, and the above bijection is an isomorphism of groups. 

In particular, the above results show that for $R$ complete and $X$
simple, the reduction mod $\mathfrak{m}$ maps the points $X(R)$
onto the points $(\mathcal{F}_{\mathfrak{m}}X)(k)$. This is a generalisation
of Hensel's lemma.

The above results can also be used to prove a generalisation of Lang's
theorem for $p$-adic group schemes. Let $R$ be the ring of integers
in a local field $F$ with finite residue field, and maximal ideal
$\mathfrak{p}$, and let $\widehat{R}^{\textrm{ur}}$ be the ring
of integers in $\widehat{F}^{\textrm{ur}}$, the completion of a maximal
unramified extension of $F$. Let $\phi\in\textrm{Gal}(\widehat{F}^{\textrm{ur}}/F)$
be the Frobenius element. Then $\phi$ restricts to an automorphism
of $\widehat{R}^{\textrm{ur}}$ over $R$. For any group scheme $G$
over $R$, we also denote by $\phi$ the induced map on $G(\widehat{R}^{\textrm{ur}})$.

\begin{prop*}
Let $G$ be a group scheme simple over $R$, whose reduction mod $\mathfrak{p}$
is connected. Then the mapping $g\mapsto g^{-1}\phi(g)$ of $G(\widehat{R}^{\mathrm{ur}})$
into itself, is surjective.
\end{prop*}
\begin{proof}
Cf. \cite{greenberg2}, \S3.
\end{proof}

\section{\label{EDL-varietiessec}Extended Deligne-Lusztig varieties}

For any discrete valuation field $F$ we denote by $\mathcal{O}_{F}$
its ring of integers, and by $\mathfrak{p}_{F}$ its maximal ideal.
For any integer $r\geq1$ we use the notation $\mathcal{O}_{F,r}$
for the quotient ring $\mathcal{O}_{F}/\mathfrak{p}_{F}^{r}$. 

Let $F$ be a local field with finite residue field $\mathbb{F}_{q}$.
We fix an algebraic closure of $F$ in which all algebraic extensions
are taken. Denote by $F^{\textrm{ur}}$ the maximal unramified extension
of $F$ with residue field $\mathbb{F}$, an algebraic closure of
$\mathbb{F}_{q}$.

Let $L_{0}$ be a finite totally ramified Galois extension of $F$,
and set $L=L_{0}^{\textrm{ur}}$, with $\Gamma=\textrm{Gal}(L/F)$.
Then $L$ is a henselian discrete valuation field with algebraically
closed residue field $\mathbb{F}$. We have the relation $\mathfrak{p}_{F}\mathcal{O}_{L}=\mathfrak{p}_{L}^{e}$,
where $e$ is the ramification index, and since $L_{0}/F$ is totally
ramified we have $e=[L_{0}:F]$. We may identify $L$ with a finite
extension of $F^{\textrm{ur}}$ (cf. \cite{ivan} ch. II, sect. 4),
and thus the residue field of $L$ is the same as that of $F^{\textrm{ur}}$.

Assume that $X$ is a scheme of finite type over $\mathcal{O}_{L}$.
In the context of the preceding section, we take $\mathcal{O}_{L}$
for the ring $R$. For any integer $r\geq1$, we define\[
X_{L,r}=\mathcal{F}_{\mathcal{O}_{L,r}}(X\underset{\ \,\mathcal{O}_{L}}{\times}\mathcal{O}_{L,r}).\]
 Conforming to the notation of \cite{lusztig-preprint-published}
and \cite{alex_paper1}, we set $X_{r}=X_{F^{\textrm{ur}},r}$. Since
$L$ and $F^{\textrm{ur}}$ have the same residue field, we have $X_{L,1}=X_{1}$. 

Let $\mathbf{G}$ be a connected reductive affine algebraic group
over $F$. Then $\mathbf{G}$ can be identified with its corresponding
affine group scheme of finite type over $F$. Extending scalars to
$L$, and using the inclusion $\mathcal{O}_{L}\rightarrow L$, we
consider $\mathbf{G}$ as a group scheme over $\mathcal{O}_{L}$.

We shall assume that $\mathbf{G}_{1}$, the reduction of $\mathbf{G}$
modulo $\mathfrak{p}_{F^{\textrm{ur}}}$, is a connected reductive
group over the residue field $\mathbb{F}$. This condition means that
for all $r>1$, we have compatibility with the case $r=1$. The condition
is satisfied in particular when $\mathbf{G}$ is a Chevalley group.

The condition implies that $\mathbf{G}$ is simple over $\mathcal{O}_{L}$,
so according to the results of the previous section, each $\mathbf{G}_{L,r}$
is an irreducible reduced affine group scheme of finite type over
$\mathbb{F}$, i.e. a connected affine algebraic group. 

Every automorphism $\sigma\in\Gamma$ stabilises $\mathcal{O}_{L}$
and $\mathfrak{p}_{L}^{r}$, respectively (cf. \cite{ivan}, chap.
II, Lemma 4.1). Therefore, each $\sigma\in\Gamma$ defines a morphism
of $\mathcal{O}_{F}$-algebras $\sigma:\mathcal{O}_{L,r}\rightarrow\mathcal{O}_{L,r}$.
Thus, by the results of the previous section, $\mathbf{G}_{L,r}$
carries a natural structure of linear algebraic group over the residue
field $\mathbb{F}$, and each $\sigma\in\Gamma$ induces a homomorphism
$\sigma:\mathbf{G}_{L,r}\rightarrow\mathbf{G}_{L,r}$ with respect
to this structure. In the following, we will use $\phi$ to denote
both the Frobenius element in $\textrm{Gal}(L/L_{0})$, and its lift
to $\Gamma$. Note that in compliance with this notation, Frobenius
morphisms on algebraic groups will in this paper always be denoted
by $\phi$.

Let $G_{r}$ denote the finite group of $\mathbb{F}_{q}$-points of
the variety $\mathcal{F}_{\mathcal{O}_{F,r}}(\mathbf{G}\times_{\mathcal{O}_{F}}\mathcal{O}_{F,r})$.
In \cite{lusztig-preprint-published} and \cite{alex_paper1}, the
group $G_{r}$ was identified with the fixed points of $\mathbf{G}_{r}$
under the Frobenius map. However, this is not the only way to realise
$G_{r}$ as a group of fixed points of an algebraic group. The following
assertion makes this more precise.

\begin{lem}
\label{lem:For-every-}For every $r\geq1$, we have $G_{r}=\mathbf{G}_{L,r'}^{\Gamma}$
if and only if $(r-1)e<r'\leq re$.
\end{lem}
\begin{proof}
For any $r'\geq1$ it is clear that $\mathcal{O}_{L,r'}^{\Gamma}=\mathcal{O}_{F,r}$,
where $r$ is the largest integer such that $\mathcal{O}_{F,r}\subseteq\mathcal{O}_{L,r'}$.
This implies that $\mathbf{G}_{L,r'}^{\Gamma}=G_{r}$, where $r$
is the largest integer such that $G_{r}\subseteq\mathbf{G}_{L,r'}$.
Now this happens exactly when $\mathfrak{p}_{F}^{r-1}\mathcal{O}_{L}\supseteq\mathfrak{p}_{L}^{r'-1}\nsubseteq\mathfrak{p}_{F}^{r}\mathcal{O}_{L}$,
i.e. when $(r-1)e\leq r'-1<re$, or equivalently $(r-1)e<r'\leq re$.
\end{proof}
Let $\Delta$ be a subset of $\Gamma$, and denote by $\{ Y_{\sigma}\mid\sigma\in\Delta\}$
a family of locally closed subsets of $\mathbf{G}_{L,r}$. For each
such family, we can define a variety\[
X=\{ g\in\mathbf{G}_{L,r}\mid g^{-1}\sigma(g)\in Y_{\sigma},\textrm{ }\forall\,\sigma\in\Delta\}.\]
Let $\left\langle \Delta\right\rangle $ denote the subgroup of $\Gamma$,
generated by $\Delta$. Then the group $\mathbf{G}_{L,r}^{\Delta}=\mathbf{G}_{L,r}^{\langle\Delta\rangle}$
clearly acts on $X$ by left multiplication. If each $Y_{\sigma}$
is normalised by a subgroup $T_{\sigma}$ of $\mathbf{G}_{L,r}$,
then there is an action of $\bigcap_{\sigma\in\Delta}(T_{\sigma})^{\sigma}$
on the variety, by right multiplication. If each $Y_{\sigma}$ is
stable under both left and right multiplication by a subgroup $U_{\sigma}$
of $\mathbf{G}_{L,r}$, and each $U_{\sigma}$ is $\sigma$-stable,
then clearly $\bigcap_{\sigma\in\Delta}U_{\sigma}$ acts on $X$ by
right multiplication. Moreover, if each $Y_{\sigma}$ is $\left\langle \Delta\right\rangle $-stable,
then we have an action of $\left\langle \Delta\right\rangle $ on
the variety.

Note that the varieties $\tilde{X}(\dot{w})$ of Deligne and Lusztig
(cf. \cite{delignelusztig}) appear as special cases of the above
construction. Namely, they are given by the specifications $r=1$,
$L=F^{\textrm{ur}}$, $\Delta=\{\phi\}$, and $Y_{\phi}=\dot{w}\mathbf{U}$,
where $\mathbf{U}$ is the unipotent radical of a Borel subgroup,
and $\dot{w}$ is a lift of an element $w$ in the Weyl group. Under
similar assumptions, but with $r\geq1$ and $Y_{\phi}=\dot{w}\mathbf{U}_{r}$,
we obtain the varieties considered by Lusztig in \cite{lusztig},
\S4 and \cite{lusztig-preprint-published}.

\begin{rem*}
We call the varieties defined in this section ``extended Deligne-Lusztig
varieties'', both because they correspond to an extension of the
maximal unramified extension, and because there are at least three
other generalisations of (certain) Deligne-Lusztig varieties, neither
of which is in the direction given here. One of these is the varieties
of Deligne associated to elements in certain braid monoids (cf. \cite{deligne-braid});
another is the affine Deligne-Lusztig varieties of Kottwitz and Rapoport
(cf.~\cite{rapoport_affine_DL}), and a third is the varieties of
Digne and Michel \cite{non_connected_DL}, defined with respect to
not necessarily connected, reductive groups. 
\end{rem*}

\section{\label{sec:An-example}An example}

Let $\mathbf{G}=\textrm{SL}_{2}$ and suppose that $F$ is of positive
characteristic with $q$ odd. Let $\epsilon$ denote a prime element
in $F$, and take $L_{0}=F[\sqrt{\epsilon}]$. Then $\Gamma$ is topologically
generated by two elements: the Frobenius automorphism $\phi$, and
an involution $\sigma$, given by $\sigma(a_{0}+a_{1}\sqrt{\epsilon})=a_{0}-a_{1}\sqrt{\epsilon}$. 

By Lemma \ref{lem:For-every-}, the smallest value of $r$ for which
$\mathbf{G}_{L,r}^{\Gamma}=G_{2}$, is $r=3$. Thus from now on, we
assume that $r=3$. We let $\Delta=\{\phi,\sigma\}$, and specify
the family $\{ Y_{\phi},Y_{\sigma}\}$ so that\[
Y_{\phi}=\begin{pmatrix}1 & 0\\
\epsilon & 1\end{pmatrix}\mathbf{U}_{L,3},\qquad Y_{\sigma}=\begin{pmatrix}1 & 0\\
\sqrt{\epsilon} & 1\end{pmatrix}\mathbf{U}_{L,3}.\]
Note that since $\mathbf{U}_{L,3}$ is closed, the same holds for
the translation by any element in $\mathbf{G}_{L,3}$. With the above
specifications, the resulting extended Deligne-Lusztig variety is\[
X_{L}:=\big\{ g\in\mathbf{G}_{L,3}\mid g^{-1}\phi(g)\in\big(\begin{smallmatrix}1 & 0\\
\epsilon & 1\end{smallmatrix}\big)\mathbf{U}_{L,3},\ \  g^{-1}\sigma(g)\in\big(\begin{smallmatrix}1 & 0\\
\sqrt{\epsilon} & 1\end{smallmatrix}\big)\mathbf{U}_{L,3}\big\}.\]
Note that this variety carries a left action of the group $G_{2}$,
and right actions of the groups \[
\mathbf{U}_{2}^{1}=\Big\{\!\!\begin{pmatrix}1 & x\epsilon\\
0 & 1\end{pmatrix}\mid x\in\mathbb{F}\Big\},\quad A=\Big\{\!\!\begin{pmatrix}\pm1+a\epsilon & 0\\
0 & \pm1-a\epsilon\end{pmatrix}\mid a\in\mathbb{F}_{q}\Big\}.\]
We wish to describe the variety $X$ more explicitly. If we let\[
g=\begin{pmatrix}a_{0}+a_{1}\sqrt{\epsilon}+a_{2}\epsilon & b_{0}+b_{1}\sqrt{\epsilon}+b_{2}\epsilon\\
c_{0}+c_{1}\sqrt{\epsilon}+c_{2}\epsilon & d_{0}+d_{1}\sqrt{\epsilon}+d_{2}\epsilon\end{pmatrix}\in\mathbf{G}_{L,3},\]
then $g^{-1}\phi(g)=\begin{pmatrix}a_{11} & a_{12}\\
a_{21} & a_{22}\end{pmatrix}$, where\begin{eqnarray*}
a_{11} & = & d_{0}a_{0}^{q}-b_{0}c_{0}^{q}+(-b_{0}c_{1}^{q}-b_{1}c_{0}^{q}+d_{0}a_{1}^{q}+d_{1}a_{0}^{q})\sqrt{\epsilon}+\\
 &  & (-b_{0}c_{2}^{q}-b_{1}c_{1}^{q}-b_{2}c_{0}^{q}+d_{0}a_{2}^{q}+d_{1}a_{1}^{q}+d_{2}a_{0}^{q})\epsilon,\\
a_{12} & = & d_{0}b_{0}^{q}-b_{0}d_{0}^{q}+(-b_{0}d_{1}^{q}-b_{1}d_{0}^{q}+d_{0}b_{1}^{q}+d_{1}b_{0}^{q})\sqrt{\epsilon}+\\
 &  & (-b_{0}d_{2}^{q}-b_{1}d_{1}^{q}-b_{2}d_{0}^{q}+d_{0}b_{2}^{q}+d_{1}b_{1}^{q}+d_{2}b_{0}^{q})\epsilon,\\
a_{21} & = & -c_{0}a_{0}^{q}+a_{0}c_{0}^{q}+(a_{0}c_{1}^{q}+a_{1}c_{0}^{q}-c_{0}a_{1}^{q}-c_{1}a_{0}^{q})\sqrt{\epsilon}+\\
 &  & (a_{0}c_{2}^{q}+a_{1}c_{1}^{q}+a_{2}c_{0}^{q}-c_{0}a_{2}^{q}-c_{1}a_{1}^{q}-c_{2}a_{0}^{q})\epsilon,\\
a_{22} & = & -c_{0}b_{0}^{q}+a_{0}d_{0}^{q}+(a_{0}d_{1}^{q}+a_{1}d_{0}^{q}-c_{0}b_{1}^{q}-c_{1}b_{0}^{q})\sqrt{\epsilon}+\\
 &  & (a_{0}d_{2}^{q}+a_{1}d_{1}^{q}+a_{2}d_{0}^{q}-c_{0}b_{2}^{q}-c_{1}b_{1}^{q}-c_{2}b_{0}^{q})\epsilon.\end{eqnarray*}
Similarly, $g^{-1}\sigma(g)=\begin{pmatrix}b_{11} & b_{12}\\
b_{21} & b_{22}\end{pmatrix}$, where\begin{eqnarray*}
b_{11} & = & 1+(b_{0}c_{1}-b_{1}c_{0}-d_{0}a_{1}+d_{1}a_{0})\sqrt{\epsilon}+\\
 &  & (-b_{0}c_{2}+b_{1}c_{1}-b_{2}c_{0}+d_{0}a_{2}-d_{1}a_{1}+d_{2}a_{0})\epsilon,\\
b_{12} & = & 2(d_{1}b_{0}-d_{0}b_{1})\sqrt{\epsilon},\\
b_{21} & = & 2(c_{0}a_{1}-c_{1}a_{0})\sqrt{\epsilon},\\
b_{22} & = & 1+(-d_{1}a_{0}+d_{0}a_{1}+b_{1}c_{0}-b_{0}c_{1})\sqrt{\epsilon}+\\
 &  & (-b_{0}c_{2}+b_{1}c_{1}-b_{2}c_{0}+d_{0}a_{2}-d_{1}a_{1}+d_{2}a_{0})\epsilon.\end{eqnarray*}
Hence, the condition for $g\in X$ becomes\[
\begin{cases}
d_{0}a_{0}^{q}-b_{0}c_{0}^{q}+(-b_{0}c_{1}^{q}-b_{1}c_{0}^{q}+d_{0}a_{1}^{q}+d_{1}a_{0}^{q})\sqrt{\epsilon}+\\
(-b_{0}c_{2}^{q}-b_{1}c_{1}^{q}-b_{2}c_{0}^{q}+d_{0}a_{2}^{q}+d_{1}a_{1}^{q}+d_{2}a_{0}^{q})\epsilon=1,\\
-c_{0}a_{0}^{q}+a_{0}c_{0}^{q}+(a_{0}c_{1}^{q}+a_{1}c_{0}^{q}-c_{0}a_{1}^{q}-c_{1}a_{0}^{q})\sqrt{\epsilon}+\\
(a_{0}c_{2}^{q}+a_{1}c_{1}^{q}+a_{2}c_{0}^{q}-c_{0}a_{2}^{q}-c_{1}a_{1}^{q}-c_{2}a_{0}^{q})\epsilon=\epsilon,\\
1+(b_{0}c_{1}-b_{1}c_{0}-d_{0}a_{1}+d_{1}a_{0})\sqrt{\epsilon}+\\
(-b_{0}c_{2}+b_{1}c_{1}-b_{2}c_{0}+d_{0}a_{2}-d_{1}a_{1}+d_{2}a_{0})\epsilon=1,\\
2(c_{0}a_{1}-c_{1}a_{0})\sqrt{\epsilon}=\sqrt{\epsilon},\\
\det(g)=1.\end{cases}\]
Note that we have omitted two redundant equations. Now, the above
system is equivalent to\[
\begin{cases}
d_{0}a_{0}^{q}-b_{0}c_{0}^{q}=1,\  a_{0}c_{0}^{q}=c_{0}a_{0}^{q},\  a_{0}d_{0}-b_{0}c_{0}=1, & (1)\\
d_{0}a_{1}^{q}+d_{1}a_{0}^{q}=b_{0}c_{1}^{q}+b_{1}c_{0}^{q},\  d_{0}a_{2}^{q}+d_{1}a_{1}^{q}+d_{2}a_{0}^{q}=b_{0}c_{2}^{q}+b_{1}c_{1}^{q}+b_{2}c_{0}^{q},\ \  & (2)\\
a_{0}c_{1}^{q}+a_{1}c_{0}^{q}=c_{0}a_{1}^{q}+c_{1}a_{0}^{q},\  a_{0}c_{2}^{q}+a_{1}c_{1}^{q}+a_{2}c_{0}^{q}-c_{0}a_{2}^{q}-c_{1}a_{1}^{q}-c_{2}a_{0}^{q}=1,\\
b_{0}c_{1}+d_{1}a_{0}=b_{1}c_{0}+d_{0}a_{1},\  d_{0}a_{2}-d_{1}a_{1}+d_{2}a_{0}=b_{0}c_{2}-b_{1}c_{1}+b_{2}c_{0}, & (4)\\
2(c_{0}a_{1}-c_{1}a_{0})=1,\\
a_{1}d_{0}+a_{0}d_{1}=b_{0}c_{1}+b_{1}c_{0},\  a_{0}d_{2}+a_{1}d_{1}+a_{2}d_{0}=b_{0}c_{2}+b_{1}c_{1}+b_{2}c_{0}. & (6)\end{cases}\]
From (1), it is easy to deduce that $a_{0}^{q}=a_{0},\  c_{0}^{q}=c_{0}$.
Using this, and subtracting the first equation (6) from the first
equations in (2) and (4), and respectively, the second equation (6),
from the second equations in (2) and (4), yields the equivalent system\[
\begin{cases}
a_{0}^{q}=a_{0},\  c_{0}^{q}=c_{0},\  a_{0}d_{0}-b_{0}c_{0}=1,\\
d_{0}(a_{1}^{q}-a_{1})=b_{0}(c_{1}^{q}-c_{1}),\  d_{0}(a_{2}^{q}-a_{2})+d_{1}a_{1}^{q}-a_{1}d_{1}=b_{0}(c_{2}^{q}-c_{2})+b_{1}c_{1}^{q}-b_{1}c_{1},\  & (2)\\
a_{0}(c_{1}^{q}-c_{1})=c_{0}(a_{1}^{q}-a_{1}),\  a_{0}(c_{2}^{q}-c_{2})-c_{0}(a_{2}^{q}-a_{2})+a_{1}c_{1}^{q}-c_{1}a_{1}^{q}=1, & (3)\\
a_{0}d_{1}=b_{1}c_{0},\  a_{1}d_{1}=b_{1}c_{1},\\
2(c_{0}a_{1}-c_{1}a_{0})=1,\\
a_{1}d_{0}+a_{0}d_{1}=b_{0}c_{1}+b_{1}c_{0},\  a_{0}d_{2}+a_{1}d_{1}+a_{2}d_{0}=b_{0}c_{2}+b_{1}c_{1}+b_{2}c_{0}.\end{cases}\]
Now, the first two equations in (2) and (3), together with $a_{0}d_{0}-b_{0}c_{0}=1$,
imply that $a_{1}^{q}=a_{1},\  c_{1}^{q}=c_{1}$. This simplifies
the other equations, so that we get\[
\begin{cases}
a_{0}^{q}=a_{0},\  c_{0}^{q}=c_{0},\  a_{0}d_{0}-b_{0}c_{0}=1,\\
a_{1}^{q}=a_{1},\  c_{1}^{q}=c_{1},\ \\
d_{0}(a_{2}^{q}-a_{2})=b_{0}(c_{2}^{q}-c_{2}),\  a_{0}(c_{2}^{q}-c_{2})-c_{0}(a_{2}^{q}-a_{2})=1,\\
a_{0}d_{1}=b_{1}c_{0},\  a_{1}d_{1}=b_{1}c_{1},\\
2(c_{0}a_{1}-c_{1}a_{0})=1,\\
a_{1}d_{0}+a_{0}d_{1}=b_{0}c_{1}+b_{1}c_{0},\  a_{0}d_{2}+a_{1}d_{1}+a_{2}d_{0}=b_{0}c_{2}+b_{1}c_{1}+b_{2}c_{0}.\end{cases}\]
Similarly, rewriting the equations in the third row, we get the equivalent
system\[
\begin{cases}
a_{0}^{q}=a_{0},\  c_{0}^{q}=c_{0},\  a_{0}d_{0}-b_{0}c_{0}=1,\\
a_{1}^{q}=a_{1},\  c_{1}^{q}=c_{1},\ \\
a_{2}^{q}-a_{2}=b_{0},\  c_{2}^{q}-c_{2}=d_{0},\\
a_{0}d_{1}=b_{1}c_{0},\  a_{1}d_{1}=b_{1}c_{1}, & (4)\\
2(c_{0}a_{1}-c_{1}a_{0})=1, & (5)\\
a_{1}d_{0}+a_{0}d_{1}=b_{0}c_{1}+b_{1}c_{0},\  a_{0}d_{2}+a_{1}d_{1}+a_{2}d_{0}=b_{0}c_{2}+b_{1}c_{1}+b_{2}c_{0}.\end{cases}\]
Now, the equations in (4) yield $a_{0}d_{1}=b_{1}c_{0}\Rightarrow a_{0}d_{1}a_{1}=a_{1}b_{1}c_{0}\Rightarrow a_{0}b_{1}c_{1}=a_{1}b_{1}c_{0}\Rightarrow b_{1}(a_{0}c_{1}-a_{1}c_{0})=0$,
and so by (5), we have $b_{1}=0$. Similarly, $d_{1}=0$.

Hence, our system of equations is equivalent to\[
\begin{cases}
a_{0}^{q}=a_{0},\  c_{0}^{q}=c_{0},\  a_{0}d_{0}-b_{0}c_{0}=1,\\
a_{1}^{q}=a_{1},\  c_{1}^{q}=c_{1},\  b_{1}=0,\  d_{1}=0,\\
a_{2}^{q}-a_{2}=b_{0},\  c_{2}^{q}-c_{2}=d_{0},\\
2(c_{0}a_{1}-c_{1}a_{0})=1,\\
a_{1}d_{0}=b_{0}c_{1},\  a_{0}d_{2}+a_{2}d_{0}=b_{0}c_{2}+b_{2}c_{0}.\end{cases}\]

Now consider the action of the group $\mathbf{U}_{2}^{1}$. If $u=\begin{pmatrix}1 & x\epsilon\\
0 & 1\end{pmatrix}\in\mathbf{U}_{2}^{1}$, then \[
gu=\begin{pmatrix}a_{0}+a_{1}\sqrt{\epsilon}+a_{2}\epsilon & b_{0}+b_{1}\sqrt{\epsilon}+(b_{2}+a_{0}x)\epsilon\\
c_{0}+c_{1}\sqrt{\epsilon}+c_{2}\epsilon & d_{0}+d_{1}\sqrt{\epsilon}+(d_{2}+c_{0}x)\epsilon\end{pmatrix}.\]
Thus the set of orbits $X_{L}/\mathbf{U}_{2}^{1}$ can be identified
with the set of points \\
$(a_{0},b_{0},c_{0},d_{0},a_{1},c_{1},a_{2},c_{2})\in\mathbb{F}^{8}$,
such that\[
\begin{cases}
a_{0}^{q}=a_{0},\  c_{0}^{q}=c_{0},\  a_{0}d_{0}-b_{0}c_{0}=1,\\
a_{1}^{q}=a_{1},\  c_{1}^{q}=c_{1},\ \\
a_{2}^{q}-a_{2}=b_{0},\  c_{2}^{q}-c_{2}=d_{0}, & (3)\\
2(c_{0}a_{1}-c_{1}a_{0})=1,\\
a_{1}d_{0}=b_{0}c_{1}.\end{cases}\]
Using the equations in (3), we eliminate $b_{0}$ and $d_{0}$, and
we can identify the variety with the set of points $(a_{0},c_{0},a_{1},c_{1},a_{2},c_{2})\in\mathbb{F}^{6}$
such that\[
\begin{cases}
a_{0}^{q}=a_{0},\  c_{0}^{q}=c_{0},\ \\
a_{1}^{q}=a_{1},\  c_{1}^{q}=c_{1},\ \\
2(c_{0}a_{1}-c_{1}a_{0})=1,\\
a_{0}(c_{2}^{q}-c_{2})-c_{0}(a_{2}^{q}-a_{2})=1,\  a_{1}(c_{2}^{q}-c_{2})=c_{1}(a_{2}^{q}-a_{2}). & (4)\end{cases}\]
The equations (4) can be rewritten so as to give the following system\[
\begin{cases}
a_{0}^{q}=a_{0},\  c_{0}^{q}=c_{0},\ \\
a_{1}^{q}=a_{1},\  c_{1}^{q}=c_{1},\ \\
2(c_{0}a_{1}-c_{1}a_{0})=1,\\
f_{1}=a_{0}c_{2}-c_{0}a_{2},\  f_{2}=a_{1}c_{2}-c_{1}a_{2},\\
f_{1}^{q}-f_{1}=1,\  f_{2}^{q}=f_{2}.\end{cases}\]
In what follows, we will denote by $Y_{L}$ the affine variety defined
by the above system of equations. To recap, we have shown the following

\begin{prop}
The quotient $X_{L}/\mathbf{U}_{2}^{1}$ is canonically isomorphic
to the affine variety $Y_{L}$.
\end{prop}
Since $\mathbf{U}_{2}^{1}$ is itself isomorphic to the affine space
$\mathbb{A}^{1}$, the irreducible representations realised in the
cohomology of $Y_{L}$ are the same as those realised by $X_{L}$
(cf. \cite{alex_paper1}, Lemma 3.1).

Following \cite{alex_paper1}, we use $(\mathbf{S},\mathbf{S})$ to
denote the subgroup of $\mathbf{G}_{2}$ consisting of matrices of
the form\[
\begin{pmatrix}1+x\epsilon & y\epsilon\\
0 & 1-x\epsilon\end{pmatrix}.\]
 Recall that the group $A=\Big\{\!\!\begin{pmatrix}\pm1+a\epsilon & 0\\
0 & \pm1-a\epsilon\end{pmatrix}\mid a\in\mathbb{F}_{q}\Big\}$ acts on $X_{L}$ by right translations.

We will now show, using the results of \cite{alex_paper1}, that the
extended Deligne-Lusztig variety $X_{L}$ realises all irreducible
representations of $G_{2}$ of dimension $(q^{2}-1)/2$. More precisely,
we show

\begin{thm}
\label{th:main}Let $Y=\{ g\in\mathbf{G}_{2}\mid g^{-1}F(g)\in(\mathbf{S},\mathbf{S})\}$.
Then there is an isomorphism\[
\alpha:Y/(\mathbf{S},\mathbf{S})\longiso Y_{L}/A,\]
which commutes with the action of $G_{2}$ on both varieties.
\end{thm}
\begin{proof}
The condition for an element $y=\left(\begin{smallmatrix}x_{0}+x_{1}\epsilon & y_{0}+y_{1}\epsilon\\
z_{0}+z_{1}\epsilon & w_{0}+w_{1}\epsilon\end{smallmatrix}\right)\in\mathbf{G}_{2}$ to lie in $Y$ is given by the equations\[
\begin{cases}
x_{0}^{q}=x_{0},\  y_{0}^{q}=y_{0},\  z_{0}^{q}=z_{0},\  w_{0}^{q}=w_{0},\  x_{0}w_{0}-y_{0}z_{0}=1,\\
x_{0}(z_{1}^{q}-z_{1})=z_{0}(x_{1}^{q}-x_{1}),\\
x_{1}w_{0}+x_{0}w_{1}=y_{0}z_{1}+y_{1}z_{0},\end{cases}\]
which can be rewritten as\[
\begin{cases}
x_{0}^{q}=x_{0},\  y_{0}^{q}=y_{0},\  z_{0}^{q}=z_{0},\  w_{0}^{q}=w_{0},\  x_{0}w_{0}-y_{0}z_{0}=1,\\
f=x_{0}z_{1}-z_{0}x_{1},\  f^{q}=f,\\
x_{1}w_{0}+x_{0}w_{1}=y_{0}z_{1}+y_{1}z_{0}.\end{cases}\]
The action on $Y$ by an element $s=\left(\begin{smallmatrix}1+t\epsilon & u\epsilon\\
0 & 1-t\epsilon\end{smallmatrix}\right)\in(\mathbf{S},\mathbf{S})$ is given by\[
ys=\begin{pmatrix}x_{0}+(x_{1}+x_{0}t)\epsilon & y_{0}+(y_{1}-y_{0}t+x_{1}u)\epsilon\\
z_{0}+(z_{1}+z_{0}t)\epsilon & w_{0}+(w_{1}-w_{0}t+z_{0}u)\epsilon\end{pmatrix}.\]
Thus, the set of orbits $Y/(\mathbf{S},\mathbf{S})$ can be identified
with the set of points \\
$(x_{0},y_{0},z_{0},w_{0},f)\in\mathbb{F}^{5}$, such that\[
\begin{cases}
x_{0}^{q}=x_{0},\  y_{0}^{q}=y_{0},\  z_{0}^{q}=z_{0},\  w_{0}^{q}=w_{0},\  x_{0}w_{0}-y_{0}z_{0}=1,\\
f^{q}=f.\end{cases}\]
Now for a point $(a_{0},c_{0},a_{1},c_{1},f_{1},f_{2})\in Y_{L}$,
the action of an element$\left(\begin{smallmatrix}1+a\epsilon & 0\\
0 & 1-a\epsilon\end{smallmatrix}\right)\in A$, is given in terms of coordinates by\[
(a_{0},c_{0},a_{1},c_{1},f_{1},f_{2})\longmapsto(a_{0},c_{0},a_{1},c_{1},f_{1},f_{2}+a/2).\]
Hence, the quotient $Y_{L}/A$ can be identified with the set of points
$(a_{0},c_{0},a_{1},c_{1},f_{1})\in\mathbb{F}^{5}$ such that \[
\begin{cases}
a_{0}^{q}=a_{0},\  c_{0}^{q}=c_{0},\ \\
a_{1}^{q}=a_{1},\  c_{1}^{q}=c_{1},\ \\
2(c_{0}a_{1}-c_{1}a_{0})=1,\\
f_{1}^{q}-f_{1}=1.\end{cases}\]
Now fix an element $\xi\in\mathbb{F}$ such that $\xi^{q}-\xi=1$.
Then there is clearly an isomorphism \[
\alpha:Y/(\mathbf{S},\mathbf{S})\longrightarrow Y_{L}/A,\]
given by\[
\alpha(x_{0},y_{0},z_{0},w_{0},f)=(x_{0},y_{0},\frac{z_{0}}{2},\frac{w_{0}}{2},f+\xi).\]
Note that because of the choice of $\xi$, this isomorphism is not
canonical. It remains to show that $\alpha$ commutes with the action
of $G_{2}$ on the varieties. Thus, let $(x_{0},y_{0},z_{0},w_{0},f)\in Y/(\mathbf{S},\mathbf{S})$,
and $(a_{0},c_{0},a_{1},c_{1},f_{1})\in Y_{L}/A$. Then the action
of an element $\begin{pmatrix}g_{0}+g_{1}\epsilon & h_{0}+h_{1}\epsilon\\
i_{0}+i_{1}\epsilon & j_{0}+j_{1}\epsilon\end{pmatrix}\in G_{2}$, is given in terms of coordinates by \begin{multline*}
(x_{0},y_{0},z_{0},w_{0},f)\longmapsto(g_{0}x_{0}+h_{0}z_{0},g_{0}y_{0}+h_{0}w_{0},i_{0}x_{0}+j_{0}z_{0},i_{0}y_{0}+j_{0}w_{0},\\
f+x_{0}^{2}(g_{0}i_{1}-i_{0}g_{1})+x_{0}z_{0}(g_{0}j_{1}+h_{0}i_{1}-i_{0}h_{1}-j_{0}g_{1})+z_{0}^{2}(h_{0}j_{1}-j_{0}h_{1})),\end{multline*}
and respectively\begin{multline*}
(a_{0},c_{0},a_{1},c_{1},f_{1})\longmapsto(g_{0}a_{0}+h_{0}c_{0},g_{0}b_{0}+h_{0}d_{0},i_{0}a_{0}+j_{0}c_{0},i_{0}b_{0}+j_{0}d_{0},\\
f_{1}+a_{0}^{2}(g_{0}i_{1}-i_{0}g_{1})+a_{0}c_{0}(g_{0}j_{1}+h_{0}i_{1}-i_{0}h_{1}-j_{0}g_{1})+c_{0}^{2}(h_{0}j_{1}-j_{0}h_{1})).\end{multline*}
Thus, it is clear that the action of $G_{2}$ commutes with the isomorphism
$\alpha$, and the theorem is proved.
\end{proof}
In \cite{alex_paper1} it was shown that all irreducible representations
of $G_{2}$ of dimension $(q^{2}-1)/2$ appear in the cohomology of
the variety $Y/(\mathbf{S},\mathbf{S})$. The above theorem shows
that that this variety is $G_{2}$-isomorphic to a quotient of an
extended Deligne-Lusztig variety by a finite group. Thus the latter
variety also realises all the above representations in its cohomology.

Of course, it would be desirable to find a more conceptual proof of
Theorem~\ref{th:main} that would not make use of explicit equations
of the varieties.

\section{Towards a reformulation of Lusztig's conjecture}

It is clear that in the degree of generality of Section \ref{EDL-varietiessec},
the varieties we have defined may sometimes be empty sets. At the
other extreme, $L_{0}=F$, $\Delta=\Gamma$, and $Y_{\sigma}=\{1\}$
for all $\sigma\in\Delta$ gives a variety identical to $G_{r}$ itself,
and thus the cohomology is just the regular representation of $G_{r}$,
which is not interesting for our purposes. Thus, in order to ensure
a nontrivial theory and a suitable framework for constructing representations
of the groups $G_{r}$, it is necessary to specialise the construction.
Motivated by the construction of Deligne and Lusztig in the case $r=1$,
the results of \cite{alex_paper1}, and the example of the preceding
section, we suggest the following preliminary construction.

As before, let $L_{0}/F$ be a finite totally ramified Galois extension
of degree $e$, and $L=L_{0}^{\textrm{ur}}$ with $\Gamma=\textrm{Gal}(L/F)$.
Let $\mathbf{G}$ be a connected reductive group over $F$. According
to a result in the structure theory of reductive groups over local
fields (cf. \cite{springer-corvallis}, 4.7), there exists a Borel
subgroup $\mathbf{B}$ in $\mathbf{G}$, defined over $F^{\textrm{ur}}$.
Let $\mathbf{U}$ denote its unipotent radical; then $\mathbf{U}$
is also defined over $F^{\textrm{ur}}$. We identify $\mathbf{G}$
and $\mathbf{U}$ with their corresponding group schemes over $\mathcal{O}_{L}$,
obtained by extension of scalars.

A first naive extension of Lusztig's construction would be the following.
Fix an integer $r\geq1$, and a corresponding group $G_{r}$. Take
$\Delta$ to be any set of topological generators of $\Gamma$, and
$r'$ any integer such that $\mathbf{G}_{L,r'}^{\Delta}=G_{r}$. For
every $\sigma\in\Delta$, let $Y_{\sigma}=x_{\sigma}\mathbf{U}_{L,r'}$
for some $x_{\sigma}\in\mathbf{G}_{L,r'}$. Denote by $X_{L,r'}(\{ x_{\sigma}\})$
the extended Deligne-Lusztig variety defined by this data, with its
corresponding action of $G_{r}$. Then one can ask whether every irreducible
representation of $G_{r}$ appears in the $l$-adic cohomology of
some $X_{L,r'}(\{ x_{\sigma}\})$.

However, the cases known so far (i.e. $r=1$, and the case discussed
in Section~\ref{sec:An-example}) are consistent with a much more
specific construction. Namely, take $\Delta=\{\phi\}\cup\Delta'$,
where $\Delta'$ is a minimal set of generators of $\textrm{Gal}(L/F^{\textrm{ur}})$.
Since $\Delta$ topologically generates $\Gamma$, we have $\mathbf{G}_{L,r'}^{\Delta}=\mathbf{G}_{L,r'}^{\Gamma}$
for any $r'$. Let $r'$ be the smallest integer such that $\mathbf{G}_{L,r'}^{\Delta}=G_{r}$.
By Lemma \ref{lem:For-every-}, this means that $r'=(r-1)e+1$. Let
$x_{\phi}$ be a representative of double $\mathbf{B}_{r}$--$\mathbf{B}_{r}$-cosets
in $\mathbf{G}_{r}$, and let $x_{\sigma}$ for each $\sigma\in\Delta'$
be a representative of double $\mathbf{B}_{L,r'}$--$\mathbf{B}_{L,r'}$-cosets
in $\mathbf{G}_{L,r'}$. Given this construction, the following questions
present themselves:\medskip

\begin{quote}
Does every irreducible representation of $G_{r}$ appear in the $l$-adic
cohomology of some variety $X_{L,r'}(\{ x_{\sigma}\})$?
\end{quote}
If so:

\begin{itemize}
\item To what extent is the construction dependent of the choice of $\mathbf{B}$?
\item Is it possible to characterise precisely what kind of extensions $L_{0}/F$
that are needed? In particular, to what extent are abelian extensions
enough?
\item Is it always sufficient to take $r'=(r-1)e+1$, or do there exist
cases where we have to take some larger $r'$ such that $(r-1)e<r'\leq re$?
\item Is it always enough to take $\Delta'$ to be a minimal set of generators
of the group $\textrm{Gal}(L/F^{\textrm{ur}})$? How does the resulting
variety depend on the choice of such a $\Delta'$?
\item For each $x_{\phi}$, can the set $\{ x_{\sigma}\mid\sigma\in\Delta'\}$
be specified further?
\end{itemize}
\medskip The answer to the first question is affirmative for $r=1$,
by the work of Deligne and Lusztig \cite{delignelusztig}. It is also
affirmative for $\mathbf{G}=\textrm{SL}_{2}$, $F$ of positive characteristic,
$q$ odd, and $r=2$, by the results in \cite{lusztig-preprint-published},
sect. 3, together with the results in \cite{alex_paper1} and Theorem~\ref{th:main}
of this paper.

\bibliographystyle{alex}
\bibliography{alex}

\end{document}